\documentclass[final]{article}
\PassOptionsToPackage{disable}{todonotes}
\input{commands-2}

\title{Accuracy, Estimates, and Representation Results}

\newcommand{\score}{\mathfrak{s}}
\renewcommand{\acc}{\mathfrak{a}}

\newcommand{\entropy}{\phi}

	\newcommand{\values}{\mathsf{Values}}

	\newcommand{\ConvHull}{\mathsf{ConvHull}}
	\newcommand{\range}{\mathsf{range}}

\acknowledgement{Huge thanks to Jason Konek for asking whether the Schervish form extended to estimates which initiated this paper. Also thanks to Richard Pettigrew, Teddy Seidenfeld, Robbie Williams, for helpful remarks.}

\begin{document}
	\maketitle 
	
\begin{abstract}

Measures of accuracy usually score how accurate a specified credence depending on whether the proposition is true or false. A key requirement for such measures is strict propriety; that probabilities expect themselves to be most accurate. We discuss characterisation results for strictly proper measures of accuracy.

By making some restrictive assumptions, we present the proof of the characterisation result of \citet{schervish1989general} in an accessible way. We will also present the characterisation in terms of Bregman divergences and the relationship between the two characterisations. 

The new contribution of the paper is to show that the Schervish form characterises proper measures of accuracy for estimates of random variables more generally, by offering a converse to \citet*[Lemma 1]{schervish2014dominating}. We also provide a Bregman divergence characterisation in the estimates setting, using the close relationship between the two forms. 
\end{abstract}

\tableofcontents

	\section{Introduction}
Considerations of accuracy---the epistemic good of having credences close to truth-values---have become an important consideration in philosophy and have led to the justification of a host of epistemic norms \Citep{joyce1998nonpragmatic,joyce2009accuracy,pettigrew2016accuracy}.
For these arguments to work, the measure of accuracy which is used needs to have certain properties. The most important of these is strict propriety, which, roughly, says that probabilities expect themselves to be the most accurate. 

The mathematical setup of accuracy measures is equivalent to scoring rules or loss functions which are used to elicit someone's credence by providing a  problem in which they are incentivised to report their own credences and in machine learning. 

There are results characterising the form that such strictly proper accuracy measures must take. By making some restrictive assumptions, we present the proof of the characterisation result of \citet{schervish1989general} in an accessible way. This provides an integral form for the strictly proper accuracy measure. We will also present the characterisation in terms of Bregman divergences and the relationship between the two characterisations. 

\Cref{sec:estimates} will discuss extending these representations to the setting where we are not measuring the accuracy of a credence, which should be close to the truth-values, but instead measuring the accuracy of an estimate of a random variable more generally, which should aim to be close to the value the random variable will actually take. The main contribution of this paper is to show that the Schervish representation result can be extended to this estimates setting, by a mild extension of the result for the credences setting. That is, we provide a converse to the result of \citet[Lemma 1]{schervish2014dominating} so that we obtain a characterisation of strictly proper measures for estimates, at least given some restrictive assumptions (\cref{thm:Schervish estimates}). This result was reported in \citet{dorstetal2021deferenceDoneBetter}, and made use of to prove their Theorem 3.2.

	\section{Accuracy of a credence}\label{sec:indicators}
	How accurate is a credence value in a proposition, say $0.6$, when the proposition is true? We give a measure. 
	
	I am not here considering the accuracy of an entire credence function at a world, but just of a single proposition. 
	
	\subsection{Setup}
	\begin{setup}We give an accuracy measure to describe how accurate a credence is in a proposition when it is true/false. Formally, we have two accuracy measures, \begin{align}
	\acc_1:[0,1]\to\Re\\
	\acc_0:[0,1]\to\Re
	\end{align}
	\end{setup}

	\begin{remark}[Infinite accuracy]
	One can often allow infinite values at end-points. 
	In particular, one can allow infinite inaccuracy at the maximally far-away points (this assumes that credences can only take values in $[0,1]$, if credences can take values in $\Re$, then we cannot have infinite values and keep truth-directendess --- we can always get worse.)
	See discussion about infinity, and various other assumptions and their relationships in \citet{schervish2009proper}.
	\end{remark}
	\begin{definition}
	$\acc$ is \defemph{(strictly) proper} iff for any $p\in[0,1]$, \begin{align}
	\Exp_p\acc(x)&:=p\acc_1(x)+(1-p)\acc_0(x)
	\end{align} obtains a (unique) maximum at $x=p$.
	\end{definition}
	\begin{definition}
	$\acc$ is \defemph{(strictly) truth-directed} iff
	If $v<x<y$ or $y<x<v$ then $\acc_v(x)>\acc_v(y)$
	\end{definition}	
	\begin{proposition}\label{thm:truthdirected}
	(Strict) propriety entails (strict) truth-directedness.
	\end{proposition}
	This is \citet[Lemma A1]{schervish1989general}. I include a proof in \cref{sec:appendix trd}. I leave this outside the main body of the paper because truth directedness is incredibly plausible, and certainly more plausible than propriety as a constraint on measurements of accuracy. (Note that this is different when one is interested in elicitation directly rather than, as philosophers usually are, measurements of the epistemic value of credences.)

	\begin{remark}
	Sometimes it would be nicer to think directly about a loss function, $\score_v$, with
		\begin{align}
		\score_v(x)&=\acc_v(v)-\acc_v(x)%
		\end{align}	

		$\score_v(x)$ measures the difference between the accuracy of perfection and the accuracy of the given credence. 		
		Note that this picture requires $\score_v(v)=0$. 	
		
		This is sometimes called a ``scoring rule'', although that terminology is also simply used for inaccuracy, i.e., negative accuracy.

		By giving a strictly proper measure $\score_v$, one can arbitrarily choose values for self-accuracy $\acc_v(v)$ to obtain a strictly proper accuracy measure by 
		\begin{align}
		\acc_v(x)&=\acc_v(v)-\score_v(x)
		\end{align} 
		
		The representations are actually really directly characterising $\score$. We can talk about strict propriety etc directly of $\score$.  This is actually more commonly done in the literature. 
		
		The literature such as \cite{pettigrew2016accuracy} works with  \emph{inaccuracy}, but I work with accuracy because it more closely ties with the philosophical presentation of trying to maximise the good of having accurate credences.

	\end{remark}

	We present two representation results for accuracy measures. 
	\subsection{Schervish}\label{sec:Schervish indicators}
	\subsubsection{Schervish form}
	The central result \citet[Theorem 4.2]{schervish1989general}
	\begin{theorem}
	$\acc$ is (strictly) proper iff
	there is some measure $\lambda$ (and values $\acc_v(v)$) such that for every $x\in[0,1]$, 
	 \begin{align}\acc_0(x)&=\acc _0(0)-\int_0^x t\;\;\lambda(\dd{t})\\
		 \acc_1(x)&=\acc _1(1)-\int_x^1 1-x\;\lambda(\diff x)
		 \end{align}
		 (for strictness, it should assign positive value to each interval)

	 \end{theorem}

			\begin{setup}\label{def:switchintegral}
			When $a>b$ define the integral $$\int_a^bf(x)\diff x=-\int_b^af(x)\diff x$$ (i.e., if it's  ``wrong-way-around'' integration limits, just take negative).
			
		Note then we can re-describe the Schervish form as: \begin{equation}
			\acc_v(x)=\acc_v(v)-\int_{x}^{v}{v-t}\;\;\lambda(\dd{t}).\label{eq:Scherv}
		\end{equation}
		i.e.,  \begin{equation}
					\score_v(x)=\int_{x}^{v}{v-t}\;\;\lambda(\dd{t}).\label{eq:Scherv}
				\end{equation} where $\score_v(x):=\acc_v(v)-\acc_v(x)$. 
		 (if $x<v$, the switching limits  and absolute value signs cancel out)
			 \end{setup}
			 
			 \begin{lemma}
			 		 A useful fact, then, is 
			 		  \begin{equation}
			 		 			\acc_v(x)=\acc_v(v)-\int_{x}^{v}{v-t}\;\;\lambda(\dd{t}).\label{eq:Scherv}
			 		 		\end{equation}
			 \end{lemma}

		\begin{remark}
		If working with inaccuracy, or the scoring rule, the signs are cleanest writing it as
		\begin{equation}
			\score_v(x)=\int_{v}^{x}{t-v}\;\;\lambda(\dd{t}).
		\end{equation}%
		\end{remark} 
			 
	\subsubsection{Any such $\acc$ is proper}	
		
		\begin{lemma}\label{thm:equivs for Schervish}
		The following are equivalent:
		\begin{enumerate}
		\item \label{itm:1}Schervish form: for $v\in\set{0,1}$ and any $x\in[0,1]$, 
		\begin{align}
		&			\acc_v(v)-\acc_v(x)=\int_{x}^{v}{v-t}\;\;\lambda(\dd{t}).&&\text{\cref{eq:Scherv}}\nonumber
		\end{align}
		\item \label{itm:2}For all $x,y\in[0,1]$, (i.e., replacing $v\in\set{0,1}$ by a general $y\in[0,1]$)
				\begin{align}
						&			\acc_v(y)-\acc_v(x)=\int_{x}^{y}{v-t}\;\;\lambda(\dd{t}).\label{eq:Schervdiff}
				\end{align}
		\item \label{itm:3}For all $x,p\in[0,1]$,
				\begin{align}
				&\Exp_p\acc(p)-\Exp_p\acc(x)=\int_x^p (p-t)\;\lambda(\dd{t})\label{eq:expintegral}
				\end{align}
		\end{enumerate}
		\end{lemma}
		\begin{proof sketch}
		These follow from quite simple manipulations. To obtain \cref{itm:1}, or \cref{itm:2} from \cref{itm:3}, note that $\Exp_v\acc(x)=\acc_v(x)$. A full proof is included in \cref{sec:equivs for Schervish}.
		\end{proof sketch}
		\begin{proposition}
		If $\acc$ has Schervish form it is (strictly) proper.
		\end{proposition}
		\begin{proof}
	Suppose $x<p$. Then for any $t\in [x,p]$, $p-t>0$, so $\int_x^p(p-t)\lambda(\diff t)>0$. 
	
	Suppose $x>p$. Then for any $t\in [x,p]$, $p-t<0$, so $\int_p^x(p-t)\lambda(\diff t)<0$. But,  $\Exp_p\acc(p)-\Exp_p\acc(x)$ switches the integral bounds, i.e., involves $\int_x^p$, which is then positive by specification of wrong-way-around integrals. 
		\end{proof}
	
	\subsubsection{Schervish's representation result}
	We prove it simply for the absolutely continuous case in order to keep the proof easy to follow. The general result holds \citep[Theorem 4.2]{schervish1989general}
	\begin{proposition}\label{thm:Schervish indicators}
	If $\acc$ is strictly proper and absolutely continuous, then there is a positive function $m$ with 
		\begin{equation}
		\acc_v(x)=\acc_v(v)-\int_x^v(v-t)m(t)\dd{t}\label{eq:schervish with mass}
		\end{equation}
	\end{proposition}
	\begin{proof}
	For absolutely continuous $\acc$, by the fundamental theorem of calculus %
	\begin{equation}
	\acc_v(v)-\acc_v(x)=\int_x^v\acc_v'(t)\dd{t}\label{eq:lemma-schervish with mass}
	\end{equation}
	
	By propriety, $\Exp_t\acc(s)$ has a maximum at $s=t$, so the derivative at this point is $0$.
	\begin{equation}
	t\acc_1'(t)+(1-t)\acc_0'(t)=0\label{eq:der}
	\end{equation}

	By manipulating \cref{eq:der} 
	\begin{equation}
	\frac{\acc_0'(t)}{-t}=\frac{\acc_1'(t)}{1-t}\label{eq:eq}
	\end{equation}
	Define the function $m$ by $m(t)=\nicefrac{\acc_0'(t)}{-t}$. So $\acc_0'(t)=-tm(t)$ and $\acc
	_1'(t)=(1-t)m(t)$. So, by replacing these in \cref{eq:lemma-schervish with mass}, we obtain \cref{eq:schervish with mass}.
	
	By \cref{thm:truthdirected}, $\acc$ is strictly truth-directed, so $\acc_0'(t)<0$ and $\acc_1'(t)>0$. Thus, $m$ is positive. 
	\end{proof}
	
	\begin{remark}
	When it is not absolutely continuous we can obtain a representation of the form: \begin{equation}
		\acc_v(x)=\acc_v(v)-\int_x^v(v-t)\;\dd{\lambda (t)}
	\end{equation}
	 we just can't push the measure $\lambda$ into a mass function. The proof of this is \citep[Theorem 4.2]{schervish1989general} and instead takes the Radon Nikodym derivatives of $\acc_0$ and $\acc_1$ relative to $\acc_1-\acc_0$. 
	 
	 Schervish also shows that the finiteness assumptions can be relaxed, 
	\end{remark}

	\subsection{Bregman divergences} \label{sec:Bregman indicators}
	Strictly proper measures can also be related to Bregman divergences \citep{savage1971elicitation,gneitingRaftery2007strictly,preddetal2009probabilistic}.
	\subsubsection{Entropy and Bregman Divergence}

	\begin{definition}
	Define the \defemph{entropy} of $\acc$ as:
	\begin{align}
	\entropy_\acc(p):=\Exp_p\acc(p)=p\acc_1(p)+(1-p)\acc_0(p)\label{eq:entropy}
	\end{align}
	\end{definition}
	
	\begin{proposition}
	If $\acc$ is proper, then $\entropy_\acc$ is convex and if it is differentiable, then:
		\begin{align}
		&\Exp_p\acc(p)-\Exp_p\acc(x)=\entropy_\acc(p)-\entropy_\acc(x)-(p-x)\entropy_\acc'(x)\label{eq:Breg}
		\end{align}
			If it is not differentiable, then we have the same form, but with $\entropy_\acc'$ as some sub-gradient. 
		
		Also for $v\in\set{0,1}$,
		\begin{align}
		&\acc_v(x)=\acc_v(v)-\br*{\entropy_\acc(v)-\entropy_\acc(x)-(v-x)\entropy_\acc'(x)}\label{eq:Breg truth}
		\end{align}
		And in fact
		\begin{align}
		&\acc_v(x)=\entropy_\acc(x)+(v-x)\entropy_\acc'(x)\label{eq:Breg truth cancel}
		\end{align} 
	\end{proposition}
	\begin{figure}[h]
	\begin{center}
	\begin{tikzpicture}[scale=3]
		\draw[thick,domain=-.8:1.6,smooth,variable=\x,blue] 
		plot ({\x},{.5*\x*\x}) node[above right] {$\Exp_p\acc(p)$};
		
		\draw[thick,red,domain=-.8:1.6,smooth,variable=\x] 
		plot ({\x},{0.1*\x -0.005}) node[below right] {$\Exp_p\acc(x)$};
		
		\filldraw[black] (0.1,0.005) circle (1pt) node[above] {$\phi(x)$};
		
		\filldraw[black] (1,0.5) circle (0.7pt) node[above left] {$\phi(p)$};
		\filldraw[black] (1,0.095) circle (0.7pt) node[below] {$\phi(x)+(p-x)\phi'(x)$};
		
		\draw[|-|,shorten >=6pt,shorten <=6pt] (1,0.5) -- (1,0.095) 
		node[midway,right] {$\phi(p)-\phi(x)-(p-x)\phi'(x)$};
	\end{tikzpicture}

	\end{center}
	\caption{Divergence diagram}	\end{figure}
	
	\begin{proof}
	By strict propriety, $\Exp_p\acc(x)<\Exp_p\acc(p)=\entropy_\acc(p)$. And 
	\begin{align}
		\Exp_p\acc(x)&=p\acc_1(x)+(1-p)\acc_0(x)
	\end{align} is a linear function of $p$ (we could name it, e.g., $f_x(p)=\Exp_p\acc(x)$). So we have a linear function entirely lying below $\entropy_\acc$ and touching it just at $p$. Therefore, $\entropy_\acc$ is convex, with $f_x(p)=\Exp_p\acc(x)$ a subtangent of it at $x$. 
	
	If $\entropy_\acc$ is differentiable at $x$, then the subtangent at $x$, which is equal to $\Exp_x\acc(p)$, is given by:
	\begin{align}
	&\Exp_p\acc(x)=\entropy_\acc(x)+(p-x)\entropy_\acc'(x)\label{eq:SavageExps}
	\end{align}
	and \cref{eq:Breg}. 
	If $\entropy_\acc$ is not differentiable, then one can take the slope of $\Exp_p\acc(x)$ and observe it is a sub-gradient of $\entropy_\acc$ by propriety; that will play the role of $\entropy_\acc'$. 

\Cref{eq:Breg truth} follows immediately, putting $p\in\set{0,1}$ and observing that $\Exp_v\acc(v)=\acc_v(v)$ to get 		\begin{align}
		&\acc_v(v)-\acc_v(x)=\br*{\entropy_\acc(v)-\entropy_\acc(x)-(v-x)\entropy_\acc'(x)}\label{eq:Breg truth}
		\end{align} and then rearranging

	\end{proof}

	\begin{definition}
	A \defemph{Bregman divergence} associated with a convex function $\entropy$ is:
	\begin{align}
	\mathfrak{d}(p,x)&:=\entropy(p)-\entropy(x)-(p-x)\entropy'(x)
	\end{align}
	\end{definition}
	So this tells us that $\Exp_p\acc(p)-\Exp_p\acc(x)$ is a Bregman divergence.

	\begin{corollary}
	If $\acc$ is strictly proper, then 
	\begin{align}
		\acc_v(x)&=\entropy_\acc(x)+(v-x)\entropy_\acc'(x)
	\end{align}where $\entropy_\acc$ the entropy for $\acc$, i.e., as in \cref{eq:entropy}.
	\end{corollary}
	\begin{proof}
	$\acc_v(x)=\Exp_v\acc(x)$. And from \cref{eq:Breg}, using the fact that $\entropy(v)=\Exp_v\acc(v)$
			\begin{align}
			&\Exp_v\acc(x)=\entropy(x)+(v-x)\entropy'(x)\label{eq:Savage}
			\end{align}
	\end{proof}

	\begin{remark}
	There is an alternative proof that goes directly via rearrangments of \cref{eq:der} using the definition of entropy, but that proof doesn't directly show that it is convex. 
	\end{remark}
	
	We also have the converse, 
	\begin{proposition}
	$\acc$ is strictly proper iff there is a convex function $\entropy$ (with values $\acc_v(v)$) where:
	\begin{align}\acc_v(x):=\acc_v(v)-\br{\entropy(v)-\entropy(x)-(v-x)\entropy'(x)}
	\end{align}
	That is, the error-score is:
	\begin{align}\score_v(x)=\entropy(v)-\entropy(x)-(v-x)\entropy'(x)
		\end{align}
	\end{proposition}

	\subsection{Relationships between Bregman divergences and the Schervish form}
	\begin{lemma}\label{thm:BregDivisIntegral thm}
	For any twice-differentiable $\entropy$, 
	\begin{align}
	\int_x^p(p-t)\entropy''(t)dt&=\entropy(p)-\entropy(x)-(p-x)\entropy'(x)\label{eq:BregDivisIntegral}
	\end{align}
	\end{lemma}
	\begin{proof}
	Integration by parts tells us that $\int_x^p u(t) v'(t)\dd{t}=\left[u(p)v(p)-u(x)v(x)\right]-\int_x^p v(t) u'(t)\dd{t}$. We apply this with $u(t):=(p-t)$, and $v(t):=\phi'(t)$, observing that $u'=-1$. So:
	{	\begin{align}
				&\int_x^p(p-t)\entropy''(t)\dd{t}\\
				&=\left[(p-p)\entropy'(p)-(p-x)\entropy'(x)\right]-\int_x^p \phi'(t)\times(-1)\dd{t}&&\text{Integration by parts}\\
				&=\int_x^p\entropy'(t)\dd{t}-(p-x)\entropy'(x)\\
				&=\entropy(p)-\entropy(x)-(p-x)\entropy'(x)
				\end{align}}	
	\end{proof}
	We can also do this with a measure rather than the mass function when $\lambda$ is a measure associated with the distribution function $\entropy'$.

	\begin{lemma}For an accuracy measure, the $m$ from Schervish and $\entropy$ the entropy, we have:
	$m(t)=\entropy''(t)$. 
	\end{lemma}
	\begin{proof}
	\begin{align}
			\entropy'(x)&=\acc_1(x)-\acc_0(x)+x\acc_1'(x)+(1-x)\acc_0'(x)&&\text{product rule}\\
			&=\acc_1(x)-\acc_0(x)&&\text{\cref{eq:der}}\label{eq:g'}
			\end{align}
	And from \cref{eq:der}, 
		\begin{equation}
		\acc_1'(x)-\acc_0'(x)=\frac{\acc_0'(x)}{-x}=m(x).\qedhere
		\end{equation}
		So $\phi''(x)=m(x)$.
	\end{proof}

	\section{Accuracy of Estimates}\label{sec:estimates}

	\subsection{Setup}

	We want to consider not only credences, which are truth-value estimates, or evaluated as good or bad with their ``closeness to the truth-value of 0/1'', but also the accuracy of one's general estimates, or forecasts, for random variables more generally.

	The random variable might, for example, be representing the utility of taking some action; or it might be the number of millimetres of rain next week. The agent will provide an estimated value. Lets say that she estimates the value to be $10$, and its true value turns out to be $30$. How accurate was her estimate? This is specified by an accuracy measure. $\acc_v:\Re\to\Re$ with $\acc_v(x)$ describing the accuracy of providing an estimated value of $x$ for $V$ when the true value is $v$. Observe that, like in the earlier setting, we are assuming that accuracy measures are finite. 
	
		\begin{setup}$\Omega$ is a non-empty set, the sample space. 
			
		A \defemph{random variable} is a function $V:\Omega\to\Re$.
		
		Probability with finite support can be given by a probability mass function $p:\Omega\to[0,1]$ with $\sum_{w\in \Omega}p(w)=1$ and where only finitely many $\omega\in\Omega$ have $p(\omega)>0$. If $\Omega$ is finite, this is just the normal notion of all probabilities. 
	\end{setup}

	\begin{definition}[Expected accuracy]
		For $p$ probabilistic with finite support, 
		\begin{align}
			\Exp_p[\acc_{V}(x)]&:=\sum_{w\in\Omega}p(\omega)\acc_{V(\omega)}(x)\\
			\Exp_p[V]&:=\sum_{w\in\Omega}p(\omega) V(\omega)
		\end{align}
	\end{definition}
	\begin{definition}[Propriety] 
		$\acc$ is \defemph{(strictly) proper for $V$} iff for any probability $p$, $\Exp_p[\acc_{V(\cdot)}(x)]$ is (uniquely) maximised at $x=\Exp_p[V]$. 
	\end{definition}
A careful specification of propriety requires spelling out what counts as a probability, e.g., finite vs countable additivity etc. We can be relaxed about this because the Schervish form entails propriety for a very wide class of probabilities, including merely finitely additive ones; and to obtain the representation we only need it to be proper for the constructed probabilities which give all their weight just to two possible values of $V$. So long as ones notion of probability contains these finitely supported probabilities and is contained by the broad notion as is spelled out in \citet{schervish2014dominating}, the results will hold.

The idea of applying scoring rules to previsions, or estimates, of random variables more generally is not novel, for example  \citet{savage1971elicitation}. The philosophical literature on accuracy has, on the whole, focused on accuracy measures as applied to credences specifically, i.e., has focused on indicator variables. 
   This characterisation of accuracy measures for estimates is, however, used in \citet{dorstetal2021deferenceDoneBetter}.

	\subsection{Schervish form for estimates}
Schervish's representation very naturally extends to consider accuracy of an estimate of any random variable. This is provided in \citet[eq 1]{schervish2014dominating}, inspired by \citet[eq 4.3]{savage1971elicitation}. It just applies the same integral form as in \cref{eq:Scherv} but allows the limits to be the true values of the variable, which may not be 0 or 1.\footnote{It was Jason Konek who suggested this form to me and asked whether the Schervish representation result extends to this setting. }

\begin{equation}
	\acc_k(x)=\acc_k(k)-\int_x^k(k-t)\;\;\lambda(\dd{t})\label{eq:Scherv estimates}
\end{equation}

\Citet[Lemma 1]{schervish2014dominating} show that all $\acc$ that have this form are strictly proper. $\lambda$ must be finite on every bounded interval, and be mutually absolutely continuous with the Lebesgue integral.

We will give the proof in the simple case where $p$ has finite support. \Citet{schervish2014dominating} do not make this assumption. Moreover, they allow $V$ to be unbounded, and in fact, they show it works for mere finite additivity expectations, defining finitely additive expectations directly as linear functionals on the extended reals that are non-negative and normalised (note: the expectation may not be determined by the underlying probability);  \citet[see][A.2]{schervish2014dominating}. They restrict the notion of propriety to those expectations, $\bb{E}$, where $\bb{E}(V)$ is finite and $\bb{E}(\acc_V(x))$ are finite for some $x$. 

\begin{proposition}[{\Citet[Lemma 1]{schervish2014dominating}}]\label{thm:schervish are str proper for estimates}
	Suppose $\acc$ has the form:
	\begin{equation}
		\acc_k(x)=\acc_k(k)-\int_x^k(k-t)\;\;\lambda(\dd{t})\label{eq:Scherv estimates}
	\end{equation}
	where $\lambda$ is finite on every bounded interval and mutually absolutely continuous with Lebesgue integral; then $\acc$ is proper. 
	
(We will here prove this for finitely supported probabilities, i.e., for any finitely supported probability $p$, $\Exp_p[\acc_V(x)]$ is maximised at $x=\Exp_p[V]$.)
	
	If $\lambda$ gives positive measure to every non-generate interval, then $\acc$ is \emph{strictly} proper, i.e., this maximum is unique.
\end{proposition}
\begin{proof}
Suppose $\acc$ is of form \cref{eq:Scherv estimates}.
	First observe that for any $k$, $x$ and $y$, \begin{equation}
		\acc _{k}(y)- \acc _{k}(x)=		\int_{x}^{y}\br{k-t}\;\;\lambda(\dd{t})\\
	\end{equation} as in the proof of $1\implies 2$ in \cref{sec:equivs for Schervish}.
	
	Let $p$ be probabilistic with finite support. 	Let $e=\Exp_p[V]$. 
	\begin{align}
		&\Exp_p[\acc(e)]-\Exp_p[\acc(x)]\\
		&=\sum_w p(w)\times (\acc _{V(w)}(e)- \acc _{V(w)}(x))\\
		&=\sum_w p(w)\times \br*{\int_{x}^{e}\br{V(w)-t}\;\;\lambda(\dd{t})}\label{eq:schervish are str proper for estimates:1}\\
		&=\int_{x}^{e}\br*{\sum_w p(w)\times \br{V(w)-t}}\;\;\lambda(\dd{t})\label{eq:schervish are str proper for estimates:2}\\
		&=\int_{x}^{e}\br*{e-t}\;\;\lambda(\dd{t})\label{eq:lemma-isproper1}
	\end{align}
	
	If $x<e$, then $e-t>0$ for all $t\in[x,e)$, and thus this integral is positive.
	
	If $x>e$, then $e-t<0$ for all $t\in(e,x]$, so $$\int_{e}^x\br*{e-t}\;\;\lambda(\dd{t})<0;$$ and thus \cref{eq:lemma-isproper1}$>0$ because the integral limits are switched, as in \cref{def:switchintegral}
	
	Thus, for any $x\neq e$, $\Exp_p[\acc(e)]-\Exp_p[\acc(x)]=\int_{x}^{e}\br*{e-t}\;\;\lambda(\dd{t})>0$.  So we know that $\Exp_p[\acc(x))]$ is maximised at $x=e$, as required. 
\end{proof}

If probabilities may not be discrete, we need to replace the sum with integrals and check that they can still be exchanged, to go from \cref{eq:schervish are str proper for estimates:1} to \cref{eq:schervish are str proper for estimates:2}; this should work so long as $V$ is bounded and the probability is countably additive. The proof given by \citet[Lemma 1]{schervish2014dominating} is in a much more general setting allowing for finite additivity and unboundedness.

We will show the converse: that any proper $\acc$ has this form.

We first will make use of a lemma
\begin{definition}
	$\acc$ is \defemph{(strictly) value-directed} iff
	If $k<x<y$ or $y<x<k$ then $\acc_k(x)>\acc_k(y)$
\end{definition}	
\begin{proposition}\label{thm:valuedirected}
	(Strict) propriety entails (strict) value-directedness.
\end{proposition}
Again we relegate the proof to the appendix because we find its fiddlyness outweighs its philosophical interest, as,  for accuracy measures, value directedness can be directly motivated.

\begin{theorem}\label{thm:Schervish estimates}
		Assume $\acc_k:\ConvHull(\values(V))\to\Re$ is absolutely continuous for each $k\in\range(V)$.
		
		If $\acc$ is (strictly) proper for $V$ (we need: for finitely supported probabilities, $p$, $\Exp_p[\acc_V(x)]$ is maximised (uniquely) at $x=\Exp_p[V]$\footnote{In fact, we just need it for the $p_{k,r,t}$ constructed; these assign all their weight to two possibilities, where $V$ takes value $k$ and $r$}); then there is a (strictly) positive function $m$ such that for every $k\in\range(V)$ and $x\in\ConvHull(\range(V))$
			\begin{equation}
			\acc_k(x)=\acc_k(k)-\int_x^k(k-t)m(t)\dd{t}
		\end{equation}
		Moreover, $$m(t)=\frac{\acc_k'(t)}{k-t}$$ for all $k\in\range(V)$.
	\end{theorem}
	
	The proof is a mild extension of the proof of \cref{thm:Schervish indicators}. 
	\begin{proof}	
	
		For absolutely continuous $\acc$, by the fundamental theorem of calculus %
		\begin{equation}
			\acc_k(k)-\acc_k(x)=\int_x^k\acc_k'(t)\dd{t}\label{eq:FTC}
		\end{equation} 
		
	For $k\in\range(V)$, define $m_k: \ConvHull(\range(V))\to\Re$ by:
	\begin{equation}
		m_k(t):=\frac{\acc_k'(t)}{k-t}.
	\end{equation}
	We then immediately get
	\begin{equation}
		\acc_k(x)=\acc_k(k)-\int_x^k(k-t)m_k(t)\dd{t}
	\end{equation}
	
		Observe that $m_k$ is (strictly) positive by (strict) value-directedness of $\acc_k$. 
		
		The important part is to show that $m_k$ does not in fact depend on the choice of $k$. That is, we need to show that $m_k=m_r$ for $k,r\in\range(V)$. For this, we use the strict propriety. 
	We first show it for $t\in\ConvHull(\set{k,r})$:
	\begin{sublemma} For $t\in\ConvHull(\set{k,r})$, we have $m_k(t)=m_r(t)$ 
	\end{sublemma}\begin{proof}
		Wlog suppose $k< r$.
	Let $\omega_k\in\Omega$ s.t.~$V(\omega_k)=k$ and $\omega_r\in\Omega$ s.t.~$V(\omega_r)=r$. Since we have assumed $t \in\ConvHull(\set{k,r})$, we con consider $p^*$ with \begin{equation}
		p^*(\omega)=\begin{cases}
			\frac{t-k}{r-k}&\omega=\omega_k\\
			\frac{r-t}{r-k}&\omega=\omega_r\\
			0&\text{otherwise}
		\end{cases}
	\end{equation}
	Observe that $p^*$ is probabilistic, i.e., these are $\geq 0$ and sum to $1$. 
	
	Now, observe that $\Exp_{p^*}[V]=t$. So by propriety, the function:
	\begin{align}
		\Exp_{p^*}\acc(x)=\frac{t-k}{r-k}\acc_r(x)+\frac{r-t}{r-k}\acc_k(x)
	\end{align}is maximised at $x=t$, so its derivative is $0$ at $t$. I.e.:
	\begin{align}
		\frac{t-k}{r-k}\acc_r'(t)+\frac{r-t}{r-k}\acc_k'(t)=0\label{eq:derexps}
	\end{align}

	By manipulating \cref{eq:derexps} we can see that:\footnote{From \cref{eq:derexps} we get:\begin{align}
			&\frac{(t-k)\acc_r'(t)+(r-t)\acc_k'(t)}{r-k}=0\\
			&(t-k)\acc_r'(t)+(r-t)\acc_k'(t)=0\\
			&(r-t)\acc_k'(t)=(k-t)\acc_r'(t)\\
			&\frac{\acc_r'(t)}{r-t}=\frac{\acc_k'(t)}{k-t}
	\end{align}}
	\begin{equation}
		\frac{\acc_r'(t)}{r-t}=\frac{\acc_k'(t)}{k-t}.\label{eq:eqexps}
	\end{equation}
	I.e., we have shown that $m_k(t)=m_r(t)$ when $t\in\ConvHull(\set{r,k})$.
	\end{proof}
	
	We still need to show $m_{k'}(t)=m_{r'}(t)$ for any $t\in\ConvHull(\range(V))$ and $k',r'\in\range(V)$. (We use $k'$ and $t'$ as we need to use the lemma for other values too). We assume wlog that $k'\leq r'$. 
	
	If $k'\leq t\leq r'$, then we can immediately apply the lemma to obtain that $m_{k'}(t)=m_{r'}(t)$. 
	
	If $t\leq k'\leq r'$, then consider also some $v\in\range(V)$ with $v\leq  t\leq k'\leq r'$. So we can apply the lemma now with $v$:
	\begin{align}
		m_{r'}(t)=m_{v}(t)&&\text{Lemma with $k=v,r=r'$}\\
		m_{k'}(t)=m_{v}(t)&&\text{Lemma with $k=v,r=k'$}\\
		\text{Thus, }m_{k'}(t)=m_{r'}(t)
	\end{align}
	Similarly, if $k'\leq r'\leq t$, we consider some $ v\in\range(V)$ with $v\geq t$ and we can similarly use the lemma to show that $m_{r'}(t)=m_{+}(t)$ and $m_{k'}(t)=m_{+}(t)$, so $m_{k'}(t)=m_{r'}(t)$. 
	
	We can thus simply put $m(t):=m_k(t)$, independent of choice of $k$. 
	\end{proof}

	\citet{schervish1989general} does not need to assume absolute continuity. 
	We conjecture that this holds in the estimates setting too. %

		\begin{corollary}Assume $\acc_k:\Re\to\Re$ is absolutely continuous for each $k\in\Re$. \\
		If $\acc$ is strictly proper for all possible random variables $V$, then there is a (strictly) positive function $m$ s.t., for every $k\in\Re$ and $x\in \Re$, 
		\begin{equation}
			\acc_k(x)=\acc_k(k)-\int_x^k(k-t)m(t)\dd{t}
		\end{equation}
		Moreover, $$m(t)=\frac{\acc_k'(t)}{k-t}$$ for all $k$.
	\end{corollary}
	\begin{proof}
		This just extends the characterisation to cover the whole of $\Re$. For any $k,x$ choose some $V$ whose maximum and minimum values cover both $k$ and $x$; then use \cref{thm:Schervish estimates}.
	\end{proof}

	\begin{corollary}
	For $p$ probabilistic with $\Exp_p[V]=e$, 
		\begin{align}
		\Exp_p\acc(e)-\Exp_p\acc(x)&=\int_x^e(e-t)m(t)\dd{t}
		\end{align}
	\end{corollary}

\begin{corollary}
	Suppose each $\acc_k$ is \emph{absolutely} continuous.
		Then $\acc$ has the form \cref{eq:Scherv estimates} iff $\acc$ is strictly proper. 
\end{corollary}
\begin{proof}
	From \cref{thm:Schervish estimates,thm:schervish are str proper for estimates}
\end{proof}

\begin{remark}[Assumptions]\label{rk:assns}
	We list various of the assumptions and their statuses:%
	\begin{itemize}
		\item Accuracy values are finite
		\begin{itemize}
			\item \Citet{schervish1989general} also applies allowing $-\infty$ accuracy, so we conjecture that the same would apply here too. 
			We need to avoid having both $+\infty$ and $-\infty$ allowed, as this causes expectations to be undefined \citep{schervish2009proper}. 
			Probably for us it's fine to allow $-\infty$ so long as it only appears at endpoints, i.e., $\acc_k(x)=-\infty\implies x\in\set{v_{\min},v_{\max}}$; which we'd get by value-directedness anyway if we were to directly assume it, which I think as accuracy measures go we'd be happy to do so. 
			\item Interestingly, though, if we want to use the same accuracy measure for all possible random variables, we will consider $\acc_k(x)$ for all $x\in\Re$, so this must be finite, as no such real-valued $x$ is an end-point for all possible variables. So the restriction can be motivated. 
			\item Consider \citet[example 1]{schervish2014dominating} showing a case where it's important that we assume the measure to be finite on any bounded interval, which is closely related to the accuracy measure only being infinity at endpoints.
		\end{itemize}
		\item Restriction to estimates evaluated being in $\ConvHull(\range(V))$
		\begin{itemize}
			\item For \cref{thm:Schervish estimates}, this restriction is essential: Strict propriety doesn't give us any control over what $\acc_k(x)$ looks like for $x\notin\ConvHull(\range(V))$, so something like the Schervish representation isn't going to be applicable. However, if we impose value-directedness in general we can still, for example, obtain accuracy-dominance results without any further control over how $\acc$ looks outside of $\ConvHull(\range(V))$
		\end{itemize}
		\item One dimensionality! We're just looking at scoring a single real-valued variable at a time. It's all one-dimensional! 
		\begin{itemize}
			\item We can push it up to finitely-many multiple variables simultaneously by just using additivity. But really it would now be natural to do this with infinitely many variables. So we're in accuracy-for-infinitely-many-propositions territory. This is exactly the sort of thing that \citet{schervish2014dominating} are considering. See also Kelley and Walsh for accuracy measures for infinitely many propositions. 
		\end{itemize}
		\item For \cref{thm:Schervish estimates} we assume that $\acc$ is absolutely continuous.
		\begin{itemize}
			\item I conjecture that the \emph{absolute} continuity is not required for \cref{thm:Schervish estimates}, just as it is not in fact required for \cref{thm:Schervish indicators} in \citet[Theorem 4.2]{schervish1989general}; of course we won't be able to push it to Schervish form with a mass function, but will need to stay in the measure setting and use Radon-Nikodyn derivatives. It becomes difficult because there are now many different $\lambda_k$ rather than just the $\lambda_0$ and $\lambda_1$ of the Schervish setting. 
			\item The simple continuity part is also probably inessential, as it is for \citet[Theorem 4.2]{schervish1989general}  because everything is anyway one-sided continuous by value-directedness. However, as \citet{schervish2009proper} point out, continuity  is essential for dominance results.
		\end{itemize}
		For \cref{thm:schervish are str proper for estimates}, \citet{schervish2014dominating} use an assumption that $\lambda$ is absolutely continuous wrt the Lebesgue measure. 
		\item Probabilities are finitely supported; can be obtained instead by assuming that $\Omega$ is finite. This isn't important:
		\begin{itemize}
			\item For the result that Schervish form entails propriety, the result holds in a very general setting \citet[Lemma 1]{schervish2014dominating}; we just present the simple proof. 
			\item For the result that propriety entails the Schervish form, we only need the fact that it is proper for these finitely supported probabilities, so being proper for a wider notion of probability would suffice. 
		\end{itemize}
		\qedhere
	\end{itemize}
	
\end{remark}

	\subsection{Bregman results}
	There is a difficulty facing the Bregman strategy which is that there is now no unique definition of entropy.\footnote{For a variable $V$ which takes values $0,0.5,1$, consider $p_1[V=1]=0.5$, $p_1[V=.5]=0$, $p_1[V=0]=0.5$, or $p_2[V=1]=0$, $p_2[V=.5]=0.5$, $p_2[V=0]=0$. 
	$\Exp_{p_1}[V]=\Exp_{p_2}[V]=0.5$. But it may be that $\Exp_{p_1}\acc(0.5)\neq\Exp_{p_2}\acc(0.5)$. }

	However, we can still get the representation by going via the Schervish representation result. 
	
	\begin{proposition}
	Assume:
	\begin{itemize}
	\item $\acc$ is continuously differentiable. What weakenings would work? Need its derivative to be integrable. 
	\item $\acc$ is absolutely continuous. NB this follows from cts diff, but perhaps not from the relevant weakenings. 
	\end{itemize}
	$\acc$ is strictly proper iff there is a strictly convex function $\entropy$ (and values $\acc_v(v)$) where:
	\begin{align}\acc_v(x)=\acc_v(v)-\left[\entropy(v)-\entropy(x)-(v-x)\entropy'(x)\right]\label{eq:BregEsts}
	\end{align}
	That is, the error-score is:
	\begin{align}\score_v(x)=\entropy(v)-\entropy(x)-(v-x)\entropy'(x)
		\end{align}
	\end{proposition}
	\begin{proof}
	We show the $\implies$ direction: Assume $\acc$ is strictly proper. 
	
	Since $\acc$ is abs cts, by \cref{thm:Schervish estimates}, we have some positive $m$ with 		
	\begin{equation}
	\acc_v(v)-\acc_v(x)=\int_x^v(v-t)m(t)\dd{t}\label{eq:lem:Scherv ests}
	\end{equation} 
	By the assumption that $\acc$ is continuously differentiable, we know that $m$ is continuous, so we can find $\phi$ where $\phi''=m$. [Weakenings?]
	
	$\phi$ is (strictly) convex since its second derivative, $m$, is (strictly) positive. 
	
	 Then, using \cref{thm:BregDivisIntegral thm}, i.e., just by integration by parts, we know that 
	\begin{align}
	\int_x^v(v-t)m(t)dt&=\entropy(v)-\entropy(x)-(v-x)\entropy'(x)\label{eq:lem:BregInt}
	\end{align}\Cref{eq:BregEsts} follows immediately from \cref{eq:lem:BregInt} and \cref{eq:lem:Scherv ests}.
	\end{proof}

The expectation form is also equivalent: 
	\begin{corollary}Also iff
	\begin{align}\Exp_p\acc(e)-\Exp_p\acc(x)={\entropy(p)-\entropy(x)-(p-x)\entropy'(x)}\label{eq:BregEsts}
	\end{align}where $e:=\Exp_p[V]$. 
	\end{corollary}
	\begin{proof}
	\begin{align}
	\Exp_p\acc(e)-\Exp_p\acc(x)&=\sum_v p[V=v] \; \br*{\acc_v(e)-\acc_v(x)}\\
	&=\sum_v p[V=v] \entropy(v)-\entropy(x)-(v-x)\entropy'(x)
	\end{align}
	\end{proof}

	\begin{remark}
		For any individual pair $k$ and $r$, we can choose some $\phi$ where $\phi(k)=\phi(r)=0$ or where $\phi(k)=\acc_k(k)$ and $\phi(r)=\acc_r(r)$; but we cannot generally choose a single $\phi$ with $\phi(v)=0$ for all $\phi$. (Since $\phi$ must be convex, there can only be two values where $\phi(v)=0$!).
	\end{remark}

	Question: how close is it to \citet{savage1971elicitation}?

	\bibcommand
\appendix
	
	\section{Propriety entails truth/value directedness}\label{sec:appendix trd}
	\subsection{Truth directedness}
		\begin{proof}[Proof of \cref{thm:truthdirected}]
	
		Take $0\leq z<y\leq 1$. We will show that $\acc_1(y)>\acc_1(z)$ and $\acc_1(y)<\acc_1(z)$.

		By strict propriety, 
		\begin{align}
		&\Exp_y\acc(y)>\Exp_y\acc(z)\\
		\text{So, }& \Exp_y[\acc(y)-\acc(z)]>0\\
		\text{So, }& y\times\br{\acc_1(y)-\acc_1(z)}+(1-y)\times\br{\acc_0(z)-\acc_0(z)}>0\label{eq:trd-1}
		\end{align}
		Let \begin{align}
		c&={\acc_1(y)-\acc_1(z)}\\
		d&={\acc_0(y)-\acc_0(z)}
		\end{align}
		So from \cref{eq:trd-1}
			\begin{equation}
			yc+(1-y)d>0\label{eq:trd-2}
			\end{equation}
		
		Similarly, by strict propriety, 
		\begin{align}
			&\Exp_z\acc(z)>\Exp_z\acc(y)\\
			\text{So, }& \Exp_z[\acc(y)-\acc(z)]<0\\
			\text{So, }& z\times\br{\acc_1(y)-\acc_1(z)}+(1-z)\br{\acc_0(z)-\acc_0(z)}\\
			&zc+(1-z)d<0&&\text{definition of $c$, $d$}\label{eq:trd-3}
			\end{align}
			
			From \cref{eq:trd-2,eq:trd-3}
			\begin{align}
			&yc+(1-y)d>zc+(1-z)d\\
			\text{So, }&(y-z)c>(y-z)d\\
			\text{Thus, }&c>d&&\text{since $y>z$}
			\end{align}
			
			Thus\begin{equation}
			c=yc+(1-y)c>yc+(1-y)d>0
			\end{equation}
			using $c>d$ for the first inequality and \cref{eq:trd-2} for the second. 
			
			Thus $c>0$. I.e., $\acc_1(y)-\acc_1(z)>0$, so $\acc_1(y)>\acc_1(z)$. 
			
			Similarly, Thus\begin{equation}
					d=yd+(1-y)d<yc+(1-y)d<0
					\end{equation}
					using $c>d$ for the first inequality and \cref{eq:trd-3} for the second. 
					
					Thus $d<0$. I.e., $\acc_0(y)-\acc_0(z)>0$, so $\acc_0(y)<\acc_0(z)$. 
		\end{proof}

			\subsection{Value directedness}\label{sec:valuedirectedproof}

				\begin{proof}[Proof of \cref{thm:valuedirected}]
			Suppose $r$ and $k$ are in the range of possible values of $V$ (with $r\neq k$). Consider $a,b$ in between $r$ and $k$, so in $[r,k]$ or $[k,r]$, and $e\in\set{r,k}$.
			
			For $x$ between $r$ and $k$ define $p_x$ assigning weight $\frac{x-r}{k-r}$ to $[V=k]$ and weight $\frac{k-x}{k-r}$ to $[V=r]$. Observe that $\Exp_{p_x}V=x$. 
			 
			Take $a,b$ between $r,k$. 
			By strict propriety, $\Exp_{p_a}\acc(b)<\Exp_{p_a}\acc(a)$ and $\Exp_{p_b}\acc(b)>\Exp_b\acc(a)$. 
			
			So $\Exp_{p_a}\acc(a)-\Exp_{p_a}\acc(b)>0>\Exp_{p_b}\acc(a)-\Exp_{p_b}\acc(b)$. I.e.,
			\begin{align}
				&\frac{a-r}{k-r}\acc_k(a)+\frac{k-a}{k-r}\acc_r(a)-\frac{a-r}{k-r}\acc_k(b)-\frac{k-a}{k-r}\acc_r(b)\\
				>\,&			\frac{b-r}{k-r}\acc_k(a)+\frac{k-b}{k-r}\acc_r(a)-\frac{b-r}{k-r}\acc_k(b)-\frac{k-b}{k-r}\acc_r(b)\\
			\end{align}
			So\begin{align}
				&\frac{a-b}{k-r}(\acc_k(a)-\acc_k(b))>\frac{a-b}{k-r}(\acc_r(a)-\acc_r(b))
			\end{align}
			Suppose $a>b>e$. Then: 
			\begin{align}
				&\frac{1}{k-r}(\acc_k(a)-\acc_k(b))>\frac{1}{k-r}(\acc_r(a)-\acc_r(b))
			\end{align}Thus
			\begin{align}
				&\acc(b,e)-\acc(a,e)\\
				&=\Exp_{p_e}\acc(b)-\Exp_{p_e}\acc(a)\\
				&=\frac{e-r}{k-r}(\acc_k(b)-\acc_k(a))+\frac{k-e}{k-r}(\acc_r(b)-\acc_r(a))\label{2eq:2'}\\
				&>\frac{b-r}{k-r}(\acc_k(b)-\acc_k(a))+\frac{k-b}{k-r}(\acc_r(b)-\acc_r(a))\label{2eq:3'}\\
				&>0
			\end{align}With \cref{2eq:2'} to \cref{2eq:3'} being because $e<b$ and there is less weight on something negative and more on something positive. 
			
			Similarly, if $a<b<e$. Then \begin{align}
				&\frac{1}{{k-r}}(\acc_k(b)-\acc_k(a))>\frac{1}{{k-r}}(\acc_r(b)-\acc_r(a))
			\end{align} So the step from \cref{2eq:2'} to \cref{2eq:3'} nonetheless holds with signs reversed. This shows value directedness whenever $a,b,e$ are between $r$ and $k$
			
			By choosing appropriate $r$, we thus show that whenever $b$ moves directly towards $k$, accuracy improves. 
		\end{proof}

			\section{Schervish equivalences}			\label{sec:equivs for Schervish}
			\begin{proof}[Proof of \cref{thm:equivs for Schervish}]

			\begin{itemize}
			\item 			\ref{itm:1}$\implies$\ref{itm:2}:
						\begin{align}
						\acc_v(y)-\acc_v(x)&=\br*{\acc_v(v)-\int_{y}^{v}{v-t}\;\;\lambda(\dd{t})}-\br*{\acc_v(v)-\int_{x}^{v}{v-t}\;\;\lambda(\dd{t})}\\
						&=\br*{\int_{x}^{v}{v-t}\;\;\lambda(\dd{t})}-\br*{\int_{y}^{v}{v-t}\;\;\lambda(\dd{t})}\\
						&=\int_{x}^{y}{v-t}\;\;\lambda(\dd{t})
						\end{align}
						
			\item			\ref{itm:2}$\implies$\ref{itm:3}:
								\begin{align}
								&\Exp_p\acc(p)-\Exp_p\acc(x)\\&=p\times\br{\acc_1(p)-\acc_1(x)}+(1-p)\times \br{\acc_0(p)-\acc_0(x)}\\
								&=p\times\br*{\int_{x}^{p}{1-t}\;\;\lambda(\dd{t})}+(1-p)\times \br*{\int_{x}^{p}{0-t}\;\;\lambda(\dd{t})}&&\text{by \cref{itm:2}}\\
								&=\int_{x}^{p}\br{p\times\br{1-t}+(1-p)\times (0-t)}\;\;\lambda(\dd{t})\\
								&=\int_{x}^{p}\br{p-t}\;\;\lambda(\dd{t})
								\end{align}
						
			\item					\ref{itm:3}$\implies$\ref{itm:1}: put $p$ as either $0$ or $1$, i.e., $v$, and simply observe that:
								\begin{align}
								\Exp_v\acc(v)=\acc_v(v)\text{ and }\Exp_v\acc(x)=\acc_v(x)
								\end{align} It then follows immediately from rearranging.
			\end{itemize}
			\end{proof}
			

\begin{thebibliography}{14}
\providecommand{\natexlab}[1]{#1}
\providecommand{\url}[1]{\texttt{#1}}
\expandafter\ifx\csname urlstyle\endcsname\relax
  \providecommand{\doi}[1]{doi: #1}\else
  \providecommand{\doi}{doi: \begingroup \urlstyle{rm}\Url}\fi

\bibitem[Briggs and Pettigrew(2016)]{briggs2016accuracy}
RA~Briggs and Richard Pettigrew.
\newblock An accuracy-dominance argument for conditionalization.
\newblock \emph{No{\^u}s}, 2016.

\bibitem[Dorst et~al.(2021)Dorst, Levinstein, Salow, Husic, and
  Fitelson]{dorstetal2021deferenceDoneBetter}
Kevin Dorst, Benjamin~A Levinstein, Bernhard Salow, Brooke~E Husic, and Branden
  Fitelson.
\newblock Deference done better.
\newblock \emph{Philosophical Perspectives}, 35\penalty0 (1):\penalty0 99--150,
  2021.

\bibitem[Gneiting and Raftery(2007)]{gneitingRaftery2007strictly}
Tilmann Gneiting and Adrian~E Raftery.
\newblock Strictly proper scoring rules, prediction, and estimation.
\newblock \emph{Journal of the American statistical Association}, 102\penalty0
  (477):\penalty0 359--378, 2007.

\bibitem[Greaves and Wallace(2006)]{greaves2006justifying}
Hilary Greaves and David Wallace.
\newblock Justifying conditionalization: Conditionalization maximizes expected
  epistemic utility.
\newblock \emph{Mind}, 115\penalty0 (459):\penalty0 607--632, 2006.

\bibitem[Joyce(1998)]{joyce1998nonpragmatic}
James~M. Joyce.
\newblock {A Nonpragmatic Vindication of Probabilism}.
\newblock \emph{Philosophy of Science}, 65:\penalty0 575--603, 1998.

\bibitem[Joyce(2009)]{joyce2009accuracy}
James~M. Joyce.
\newblock Accuracy and coherence: Prospects for an alethic epistemology of
  partial belief.
\newblock In \emph{Degrees of belief}, pages 263--297. Springer, 2009.

\bibitem[Pettigrew(2014)]{pettigrew2014accuracyindifference}
Richard Pettigrew.
\newblock Accuracy, risk, and the principle of indifference.
\newblock \emph{Philosophy and Phenomenological Research}, 92\penalty0
  (1):\penalty0 35--59, 2014.

\bibitem[Pettigrew(2016)]{pettigrew2016accuracy}
Richard Pettigrew.
\newblock \emph{Accuracy and the Laws of Credence}.
\newblock Oxford University Press, 2016.

\bibitem[Predd et~al.(2009)Predd, Seiringer, Lieb, Osherson, Poor, and
  Kulkarni]{preddetal2009probabilistic}
Joel~B Predd, Robert Seiringer, Elliott~H Lieb, Daniel~N Osherson, H~Vincent
  Poor, and Sanjeev~R Kulkarni.
\newblock Probabilistic coherence and proper scoring rules.
\newblock \emph{IEEE Transactions on Information Theory}, 55\penalty0
  (10):\penalty0 4786--4792, 2009.

\bibitem[Savage(1971)]{savage1971elicitation}
Leonard~J Savage.
\newblock Elicitation of personal probabilities and expectations.
\newblock \emph{Journal of the American Statistical Association}, 66\penalty0
  (336):\penalty0 783--801, 1971.

\bibitem[Schervish(1989)]{schervish1989general}
Mark~J Schervish.
\newblock A general method for comparing probability assessors.
\newblock \emph{The Annals of Statistics}, 17.4:\penalty0 1856--1879, 1989.

\bibitem[Schervish et~al.(2009)Schervish, Seidenfeld, and
  Kadane]{schervish2009proper}
Mark~J Schervish, Teddy Seidenfeld, and Joseph~B Kadane.
\newblock Proper scoring rules, dominated forecasts, and coherence.
\newblock \emph{Decision Analysis}, 6\penalty0 (4):\penalty0 202--221, 2009.

\bibitem[Schervish et~al.(2014)Schervish, Seidenfeld, and
  Kadane]{schervish2014dominating}
Mark~J Schervish, Teddy Seidenfeld, and Joseph~B Kadane.
\newblock Dominating countably many forecasts.
\newblock \emph{The Annals of Statistics}, 42\penalty0 (2):\penalty0 728--756,
  2014.

\bibitem[Walley(2000)]{walley2000towards}
Peter Walley.
\newblock Towards a unified theory of imprecise probability.
\newblock \emph{International Journal of Approximate Reasoning}, 24\penalty0
  (2-3):\penalty0 125--148, 2000.

\end{thebibliography}
		\end{document}